\documentclass{article}
\usepackage{spconf}

\usepackage{cite}
\usepackage{amsmath,amssymb,amsfonts}
\usepackage{algorithmic}
\usepackage{graphicx}
\usepackage{array}
\usepackage{textcomp}
\usepackage{multirow}
\usepackage[table,xcdraw]{xcolor}
\usepackage{comment}
\usepackage{leftidx}
\usepackage{algorithm2e}
\usepackage{amsthm}
\usepackage{float}
\usepackage[position=top]{subfig}
\usepackage{setspace}
\RestyleAlgo{ruled}
\newcommand\ER[1]{\textcolor{black}{#1}}

\newcommand\np[1]{\textcolor{black}{#1}}

\newtheorem{theorem}{Theorem}

\newtheorem{lemma}{Lemma}
\newcommand{\yhk}{y_{h,k}}
\newcommand{\yhkun}{y_{h,k+1}}

\newcommand{\xhk}{x_{h,k}}
\newcommand{\xhkun}{x_{h,k+1}}
\newcommand{\xHl}{x_{H,k,\ell}}
\newcommand{\xHlun}{x_{H,k,\ell+1}}
\newcommand{\xHO}{x_{H,k,0}}
\newcommand{\xHm}{x_{H,k,m}}

\newcommand{\taubarh}{\bar{\tau}_{h,k}}

\newcommand{\tauhk}{\tau_{h,k}}

\newcommand{\alphalk}{\alpha_{H,k,\ell}}

\newcommand{\ehk}{c_{h,k}}
\newcommand{\Dhk}{D_{h,k}}

\newcommand{\vHk}{v_{H,k}}
\newcommand{\A}{\mathrm{A}}
\newcommand{\W}{\mathrm{W}}
\newcommand{\I}{\mathrm{I}}
\newcommand{\RR}{\mathbb{R}}

\newcommand{\prox}{\mbox{prox}}
\newcommand{\ftn}{\scriptsize}
\newcolumntype{L}[1]{>{\raggedright\let\newline\\\arraybackslash\hspace{0pt}}m{#1}}
\newcolumntype{C}[1]{>{\centering\let\newline\\\arraybackslash\hspace{0pt}}m{#1}}
\newcolumntype{R}[1]{>{\raggedleft\let\newline\\\arraybackslash\hspace{0pt}}m{#1}}
\title{Multilevel FISTA for Image Restoration}
\name{Guillaume Lauga$^{\dagger}$, Elisa Riccietti$^{\dagger}$, Nelly Pustelnik$^{\star}$, Paulo Gon\c calves$^{\dagger}$\thanks{The authors would like to thank the GdR ISIS for funding the MOMIGS project and the ANR-19-CE48-0009 Multisc'In project. We also gratefully acknowledge the support of the Centre Blaise Pascal's IT test platform at ENS de Lyon (Lyon, France) for the computing facilities. The platform operates the SIDUS \cite{quemener2013} solution developed by Emmanuel Quemener. All experiments were performed on a machine with 8-cores Intel Core i7-10700, 4.8 GHz, 62 Go of RAM.}}
\address{$^{\dagger}$Univ Lyon, Inria, EnsL, UCBL, CNRS, LIP, UMR 5668, F-69342, Lyon Cedex 07, France \\ $^{\star}$Ens de Lyon,
		CNRS, Laboratoire de Physique, F-69342, Lyon, France}

\begin{document}
\maketitle
\begin{abstract}
This paper presents a multilevel FISTA algorithm, based on the use  of the Moreau envelope to build the correction brought by the coarse models, which is easy to compute when the explicit form of the proximal operator of the considered functions is known.
This approach is supported by strong theoretical guarantees: we prove both the rate of convergence and the convergence of the iterates to a minimum in the convex case, an important result for ill-posed problems. We evaluate our approach on image restoration problems and we show that it outperforms classical FISTA for large-scale images.
\end{abstract}

\begin{keywords}
multilevel optimization, inertial methods, image restoration, proximal methods.
\end{keywords}
\section{Introduction}
Many problems in signal and image processing involve minimising a sum of a data fidelity term $f$ and a regularization function $g$, formally:
\begin{equation}
    \min_{x \in \mathcal{H}} F(x):=f(x) + g(x)
    \label{eq:optim_classic}
\end{equation}
where $\mathcal{H}$ is a real Hilbert space ($\mathcal{H} = \mathbb{R}^N$ in the following), $f:\mathcal{H} \rightarrow (-\infty,+\infty]$ and $g:\mathcal{H} \rightarrow (-\infty,+\infty]$ belong to $\Gamma_0(\mathcal{H})$ the class of convex, lower semi-continuous, and proper functions. Moreover, $f$ is assumed to be  differentiable with gradient $L_f$-Lipschitz and $F$ is supposed to be coercive. 

In the context of restoration, we aim to recover a good quality image from \np{an} 
image $z$, degraded by a linear operator and/or noise, i.e. $z = \A\bar{x}+\epsilon$, where $\A \in \mathbb{R}^{N\times N}$ models the linear degradation operator, $\epsilon$ the noise. 
To solve this ill-posed problem, we generally consider a regularized least squares formulation, where we denote $g$ the regularization function allowing us to choose the properties that we wish to impose on the solution. 
A usual choice is to apply the $l_1$-norm on the coefficients raised by a linear transformation $\W\in \mathbb{R}^{K\times N}$ (wavelets, frames, dictionary, \ldots), thus promoting the sparsity of the solution \cite{pustelnik2016}. 
Given a regularization parameter $\lambda>0$, the associated minimization problem reads:
\begin{equation}
    \widehat{x} \in \underset{x \in \mathbb{R}^N}{\textrm{Argmin}} \frac{1}{2} \Vert \A x-z \Vert_2^2 + \lambda \Vert \W x \Vert_1.
    \label{eq:optim_classique}
\end{equation}
Many algorithms have been proposed in the literature to estimate $\widehat{x}$ (cf. \cite{Combettes_P_2011,parikh2014,chambolle2016}). They suffer from the significant increase in computational time with the dimension. Preconditioning techniques can be investigated but generally require strong assumptions for the choice of the preconditioning matrix (e.g., diagonal matrix) leading to limited gains. For the solution of large-scale problems with smooth objective function, it is possible to take advantage of the local structure of the optimisation problem (cf. VMLMB \cite{VMLM-B} or 3MG \cite{3MG}). %

In this paper, we focus on a different family of approaches, the multilevel schemes, which exploit different resolutions of the same problem.  In such methods the objective function is approximated by a sequence of functions defined %
on reduced dimensional spaces (coarser spaces).%
The descent step is thus calculated %
at \np{coarser} 
levels with minimal cost and then projected to the fine levels. 

These approaches have been mainly studied for the solution of partial differential equations (PDEs), in which $f$ and $g$ are supposed to be differentiable \cite{nash2000,calandra2020}, but recently this idea has also been exploited in \cite{parpas2016,parpas2017,lauga2022} to define multilevel forward-backward proximal algorithms applicable to problem \eqref{eq:optim_classic} in the case where $g$ is non differentiable.

In this paper we propose a variant of these methods, which we call MMFISTA for \textit{Moreau Multilevel FISTA} providing a multilevel alternative to inertial strategies such as FISTA~\cite{beck2009-1,aujol2015}. %
Our framework relies on the Moreau envelope for the definition of smooth coarse 
approximations of $g$, which can be easily constructed when the proximal operator of $g$ is known in explicit form.
Furthermore, we show under mild assumptions that the convergence guarantees of FISTA hold for MMFISTA, in particular the convergence of the iterates, an important result for ill-posed problems and, to our knowledge, never established for multilevel inertial proximal methods. %

The paper is organized as follows.
In Section~\ref{sec:notre_methode}, we recall the main principles of FISTA. Then, we describe MMFISTA. In Section~\ref{sec:convergence}, we present its convergence guarantees. Finally, in Section~\ref{sec:results}, we present numerical results to confirm the good behaviour of MMFISTA in an image restoration context. 

\vspace{-1em}
\section{Multilevel FISTA}
\label{sec:notre_methode}
\noindent\textbf{FISTA} -- Among the numerous algorithms designed to solve a minimization problem of the form \eqref{eq:optim_classic}, the most standard strategy is FISTA \cite{beck2009-1}, which relies on forward-backward iterations and extrapolation steps, such that, for every $k=0,1,\ldots$
\begin{align}
    x_{k+1} & =\prox_{\tau_k g}(y_k - \tau_k \nabla f(y_k))\label{eq:fista_steps1} \\
    y_{k+1} & = x_{k+1} + \alpha_{k}(x_{k+1}-x_{k})
    \label{eq:fista_steps2}
\end{align}
where $x_0 =y_0$ and $\alpha_{k} = \frac{t_{k}-1}{t_{k+1}}$ for all $k\geq 1$. Choosing $t_k = \left(\frac{k+a-1}{a} \right)$ where $a>2$ \cite[Definition 3.1]{aujol2015} and $\tau_k \in (0, L_{f}^{-1})$ ensures various convergence guarantees (see \cite[Theorem 3.5 and 4.1]{aujol2015}). We will denote these conditions (AD) in the following. 

\noindent\textbf{Multilevel framework} -- %
The multilevel framework exploits a hierarchy of objective functions, which are representations of $F$ at different resolutions  and alternate minimization between these objective functions (following a V cycle procedure \cite{nash2000}). Without loss of generality and for the sake of clarity, we consider the two-level case: we index by $h$ (resp. $H$) all quantities defined at the fine (resp. coarse) level. We thus define %
$F_h:=F: \mathbb{R}^{N_h}\rightarrow (-\infty,+\infty]$ the objective function at the fine level where $N_h=N$, involving $f_h:=f$ and $g_h:=g$. Its approximation at the coarse level is denoted $F_H: \mathbb{R}^{N_H}\rightarrow (-\infty,+\infty]$ where  $N_H<N_h$, {which involves} $f_H$ and $g_H$. We also define transfer information ope\-rators: a linear operator $\I_h^H: \RR^{N_h} \to \RR^{N_H}$ that sends information from the fine level to the coarse level, and conversely $\I_H^h: \RR^{N_H} \to \RR^{N_h}$ that sends information from the coarse level back to the fine level.  It is classical to choose $\I_H^h = \eta (\I_h^H)^T$, with $\eta >0$.

In a multilevel scheme,  we improve the intermediate %
iterate $ \yhk$ by
performing iterations at the coarse level: %
$\yhk$ is projected to the coarse level with $I_h^H$ \eqref{eq:step1},  %
a sequence $(\xHl)_{\ell \in \mathbb{N}}$ is defined (where $k$ represents the current iteration at the fine level and $\ell$ indexes the iterations at the coarse level) such that: ${x}_{H,k,\ell+1} = \Phi_{H,k,\ell}(\xHl)$, with $\Phi_{H,k,\ell}$ any operator such that $F_H(\xHm)\leq F_H(\xHO)$ for some $m>0$. This yields after $m$ iterations at the coarse level \eqref{eq:step2} %
\np{to a step being brought back at the fine level \eqref{eq:step3}.}
%
%
Then, the generic iteration $k$ of a multilevel method reads :
\begin{subequations}
\begin{align}
     \xHO & =\I_h^H \yhk \label{eq:step1}\\
   \xHm &= \Phi_{H,k,m-1} \circ .. \circ \Phi_{H,k,0}(\xHO)\label{eq:step2}\\ 
    \bar{y}_{h,k} &= \yhk + \mbox{$\np{\taubarh}$} \I_H^h\left(\xHm-\xHO\right) \label{eq:step3}\\
    x_{h,k+1} & =\prox_{\tau_{h,k} g_h}(\bar{y}_{h,k}  - \tau_k \nabla f_h(\bar{y}_{h,k} )) \label{eq:step4}\\
    \yhkun & = \xhkun + \alpha_{h,k}(\xhkun-x_{h,k}) \label{eq:step5}
\end{align}
\end{subequations}
By taking $\xHm = \xHO$ one recovers the standard FISTA iteration.
To ensure that the correction term %
$\xHm-\xHO$, \ER{once} projected from coarse level to fine level, provides a decrease of $F_h$, we need to do appropriate choices for :
\begin{itemize}
\item the coarse model $F_H$,
\item the minimization scheme $\Phi_{H}$.
\end{itemize}
\noindent \textbf{Coarse model $F_H$} --  %
The coarse iterations are built using the Moreau envelope of $g_h$ and of its coarse approximation $g_H$. \np{The Moreau envelope provides a natural choice to extend ideas coming from the classical smooth case \cite{calandra2020} to proximal gradient methods because of its smoothness and its expression involving the proximity operator.} 
We first recall that for $\gamma>0$ and $g$ being a convex, lower semi-continuous, and proper function of $\mathcal{H}$ in $(-\infty,+\infty]$, its Moreau envelope, denoted $\leftidx{^{\gamma}}{g}$, is the convex, continuous, real-valued function defined by
\begin{equation}
    \leftidx{^{\gamma}}{g} = \inf_{y \in \mathcal{H}} g(y) + (1/2 \gamma) \Vert \cdot -y \Vert^2,
\end{equation}
which can be expressed explicitly with $\prox_{\gamma g}$ \cite[Remark 12.24]{bauschke2017}.
The gradient of $\leftidx{^{\gamma}}{g}$ is $\gamma^{-1}$-lipschitz and such that \textcolor{black}{(Prop. 12.30 in \cite{bauschke2017})}
\begin{equation}
    \nabla(\leftidx{^{\gamma}}{g}) = \gamma^{-1}(\mbox{Id}- \prox_{\gamma g}).
    \label{gradient_moreau}
\end{equation}
\noindent The coarse model $F_H$ is  defined as 
\begin{equation}
\label{eq:FH}
 F_{H}(x_H) = f_H(x_H) + g_H(x_H) + \langle \vHk,x_H \rangle
\end{equation}
where 
\vspace{-1em}
\begin{align}
\label{eq:v_Hk}
\vHk = & I_h^H \left(\nabla f_h(\yhk) + \nabla (\leftidx{^{\gamma_h}}{g}{}_h)(\yhk)\right) \nonumber\\ & -(\nabla f_H(\xHO) + \nabla(\leftidx{^{\gamma_H}}{g}{}_H)(\xHO)).
\end{align}
The third term in \eqref{eq:FH} %
is added to enforce the first order coherence between a smoothed coarse objective function \begin{equation}\label{eq:F_Hgamma} %
F_{H,\gamma_H}(x_H) = f_H(x_H) + \leftidx{^{\gamma_H}}{g}_H(x_H) + \langle \vHk,x_H \rangle
\end{equation}
\noindent and a smoothed fine objective function $F_{h,\gamma_h}$\cite{parpas2017} near $\xHO$:
\begin{equation}\label{descent}
    \nabla F_{H,\gamma_H}(\xHO) = I_h^H \nabla F_{h,\gamma_h}(\yhk).
\end{equation}
The choice of the smoothing parameters $\gamma_h$ %
and $\gamma_H$ will be discussed in Section \ref{sec:results}.
This condition ensures that if 
$\xHm-\xHO$ is a descent direction for $F_{H,\gamma_H}$ at $\xHO$, then $I_H^h(\xHm-\xHO)$ is a descent direction for $F_{h,\gamma_h}$ as well: 
\vspace{-0.5em}
$$\langle I_H^h(\xHm-\xHO),\nabla F_{h,\gamma_h}(\xhk) \rangle \leq 0.$$
%
%
%
%
%
\noindent %
According to properties of the Moreau envelope and the principles developed in \cite{beck2012}, if $\xHm-\xHO$ is a descent direction for $F_{H,\gamma_H}$, we obtain $$F_h(\yhk + \bar{\tau}_{h,k}I_H^h(x_{H,k,m}-x_{H,k,0})) \leq F_h(\yhk) + \beta \gamma_h$$
\noindent where $\beta>0$ depends on $g_h$. This ensures that $F_h$ is decreasing up to a constant $\beta \gamma_h$ (which can be made arbitrarily small) after a use of the coarse models. %
Now we show how to enforce the decrease of $F_{H,\gamma_H}$.

\noindent\textbf{Minimization operator $\Phi_{H}$ --} 
At the coarse level we can decide to consider either the non-smooth approximation \eqref{eq:FH} of the objective function or the smoothed version \eqref{eq:F_Hgamma}. Both cases lead to a decrease in $F_{H,\gamma_H}$: %
\np{indeed, taking the Moreau envelope of $g_H$ in $F_{H}(\xHm)\leq F_{H}(\xHO)$ yields $F_{H,\gamma_H}(\xHm)\leq F_{H,\gamma_H}(\xHO)$}. The two cases are linked by the same choice of the correction term to ensure the coherence between the two levels \eqref{eq:v_Hk}. 
We consider here three different strategies :
\begin{enumerate}
    \item Gradient steps on the smoothed $F_{H,\gamma_H}$: \\ $\Phi_{H,S} = \left(\mbox{Id}-\tau_H (\nabla (f_H+\leftidx{^{\gamma_H}}{g}{}_H)+v_H)\right)$ 
    \item Proximal gradient steps on the non-smooth $F_{H}$: \\ $\Phi_{H,FB} = \prox_{\tau_H g_H} \left(\mbox{Id}-\tau_H (\nabla f_H+v_H)\right)$.
    \item FISTA steps on the non-smooth $F_H$ with the previous proximal gradient step and where $\alphalk$ follows (AD) conditions. Noted $\Phi_{H,FISTA}$ in the following. 
\end{enumerate}

\noindent \textbf{Practical considerations --}
Our algorithm is based on a simple construction of %
$F_H$ and $\vHk$, as long as the computation of the associated proximal operator has an explicit form, which is a rather reasonable assumption. Our method is sketched in Algorithm \ref{alg:MMFISTA}. The step length at both levels can be selected either by fixing a value below the threshold guaranteeing convergence, defined by the Lipschitz constants associated to the functions considered, or by a linear search. The second solution is generally more costly, but may provide faster convergence in some cases. 
To ensure the convergence of the iterates, we impose at most $p$ uses of the coarse models $F_H$ (one use corresponds to a full V-scheme cycle), which is also recommended to obtain a good computation time (cf. Section \ref{sec:results}).

\begin{algorithm}
\caption{MMFISTA}\label{alg:MMFISTA}
\KwData{\ftn $x_{h,0}$, $\epsilon_h,\gamma,m,p>0$, $t_{h,0}=1$, $a>2$, $k=0$, $r=0$}
\While{$\Vert \xhkun-\xhk\Vert> \epsilon_h$}{
  \eIf{$r< p$}{
  \small
  $r=r+1$\\
    $\xHO = y_{H,k,0} = I_h^H \yhk$\\
    $\vHk = I_h^H \nabla F_{h,\gamma_h}(\yhk)- \nabla F_{H,\gamma_H}(\xHO)$ \\
    \For{$\ell=0\dots m-1$}{
    \ftn $y_{H,k,\ell+1} = \Phi_{H,k,\ell}(\xHl)$ \\
    $\xHlun = y_{H,k,\ell+1} + \alpha_{H,k,\ell}(y_{H,k,\ell+1}-y_{H,k,\ell})$ \\
    } 
    Set $\taubarh>0$, 
    $\bar{y}_{h,k} = \yhk + \taubarh I_H^h(\xHm-\xHO)$ \\
    Set $\tauhk>0$, \\
    $\xhkun = \prox_{\tauhk g_h}(\bar{y}_{h,k} - \tauhk \nabla f_h(\bar{y}_{h,k}))$
  }{\small Set $\tauhk>0$,  \\
    $\xhkun = \prox_{\tauhk g_h}(y_{h,k} - \tauhk \nabla f_h(y_{h,k}))$
  }
  \small $t_{h,k}= \left(\frac{k+a-1}{a} \right)$, $\alpha_{h,k} = \frac{t_{h,k}-1}{t_{h,k+1}}$ \\
  \small $\yhkun = \xhkun + \alpha_{h,k} (\xhkun-\xhk)$
}
\end{algorithm}
\vspace{-1em}

\begin{table*}[ht]
        \centering
        \setlength\tabcolsep{1pt}
        \begin{tabular}{|C{20mm}|C{20mm}|C{12mm}|C{12mm}|C{12mm}|C{12mm}|C{12mm}|C{5mm}|C{12mm}|C{12mm}|C{12mm}|C{12mm}|C{12mm}|}
        \hline
        \multicolumn{2}{|c|}{\ftn \textbf{Noise $\backslash$ Blur} }& \multicolumn{5}{c|}{\ftn (a) size(blur) = $[40,40]$, $\sigma$(blur) = 7.3} & & \multicolumn{5}{c|}{\ftn  (b) size(blur) = $[88,88]$, $\sigma$(blur) = 16}\\
        \hline
        \multirow{4}{*}{$(1)$ $ \sigma=0.01$}&\multicolumn{1}{|c|}{\cellcolor[HTML]{C0C0C0}\ftn FISTA CPU time}& \ftn \cellcolor[HTML]{C0C0C0}16& \cellcolor[HTML]{C0C0C0}\ftn28&\cellcolor[HTML]{C0C0C0}\ftn 42&\cellcolor[HTML]{C0C0C0}\ftn 161 &\cellcolor[HTML]{C0C0C0}\ftn 401 & &\cellcolor[HTML]{C0C0C0} \ftn 17 &\cellcolor[HTML]{C0C0C0}\ftn 30 &\cellcolor[HTML]{C0C0C0}\ftn42&\cellcolor[HTML]{C0C0C0}\ftn 148&\cellcolor[HTML]{C0C0C0}\ftn 421\\ \cline{2-7} \cline{9-13}
        &\multicolumn{1}{|c|}{\ftn$\Phi_{H,S}$}& \ftn$-20~\textcolor{red}{\bullet}$ & \ftn$-22~\textcolor{blue}{\bullet}$& \ftn$+1~\textcolor{red}{\bullet}$  & \ftn$+1~\textcolor{red}{\bullet}$  &\ftn $-1~\textcolor{red}{\bullet}$ & & \ftn$\mathbf{-51}~\textcolor{red}{\bullet}$ & \ftn$\mathbf{-44}~\textcolor{blue}{\bullet}$& \ftn$-18~\textcolor{blue}{\bullet}$  & \ftn$\mathbf{+4}~\textcolor{red}{\bullet}$  &\ftn $\mathbf{-1}~\textcolor{red}{\bullet}$ \\ \cline{2-7} \cline{9-13}
        & \multicolumn{1}{|c|}{\ftn$\Phi_{H,FB}$}& \ftn$-19~\textcolor{red}{\bullet}$& \ftn$-19~\textcolor{blue}{\bullet}$& \ftn$+5~\textcolor{red}{\bullet}$ & \ftn$\mathbf{+2}~\textcolor{red}{\bullet}$&\ftn$\mathbf{+1}~\textcolor{red}{\bullet}$& & \ftn$-50~\textcolor{red}{\bullet}$ & \ftn$-42~\textcolor{blue}{\bullet}$& \ftn$-15~\textcolor{blue}{\bullet}$   & \ftn$+6~\textcolor{red}{\bullet}$& \ftn$+1~\textcolor{red}{\bullet}$\\ \cline{2-7} \cline{9-13}
        & \multicolumn{1}{|c|}{\ftn$\Phi_{H,FISTA}$}&\ftn$\mathbf{-51}~\textcolor{red}{\bullet}$ &\ftn$\mathbf{-32}~\textcolor{blue}{\bullet}$& \ftn$\mathbf{-4}~\textcolor{blue}{\bullet}$&\ftn$+2~\textcolor{red}{\bullet}$& \ftn$+1~\textcolor{red}{\bullet}$ & & \ftn$-50~\textcolor{red}{\bullet}$&\ftn$-42~\textcolor{blue}{\bullet}$& \ftn$\mathbf{-35}~\textcolor{blue}{\bullet}$&\ftn$+8~\textcolor{red}{\bullet}$& \ftn$+1~\textcolor{red}{\bullet}$ \\ \cline{2-7} \cline{9-13}
        \hline
        \multirow{4}{*}{$(2)$ $\sigma=0.04$}&\multicolumn{1}{|c|}{\cellcolor[HTML]{C0C0C0}\ftn FISTA CPU time}&\cellcolor[HTML]{C0C0C0}\ftn 14 &\cellcolor[HTML]{C0C0C0}\ftn 22&\cellcolor[HTML]{C0C0C0}\ftn 34  &\cellcolor[HTML]{C0C0C0}\ftn 108 &\cellcolor[HTML]{C0C0C0}\ftn 220 & &\cellcolor[HTML]{C0C0C0} \ftn 15&\cellcolor[HTML]{C0C0C0}\ftn 25&\cellcolor[HTML]{C0C0C0}\ftn 34&\cellcolor[HTML]{C0C0C0}\ftn 122&\cellcolor[HTML]{C0C0C0}\ftn 315 \\ \cline{2-7} \cline{9-13}
        & \multicolumn{1}{|c|}{\ftn$\Phi_{H,S}$}& \ftn$\mathbf{-22}~\textcolor{red}{\bullet}$ & \ftn$-10~\textcolor{blue}{\bullet}$& \ftn$-1~\textcolor{red}{\bullet}$  & \ftn$\mathbf{-1}~\textcolor{red}{\bullet}$  &\ftn $-1~\textcolor{red}{\bullet}$  & & \ftn$-29~\textcolor{red}{\bullet}$ & \ftn$-25~\textcolor{blue}{\bullet}$& \ftn$-18~\textcolor{blue}{\bullet}$ & \ftn$\mathbf{+3}~\textcolor{red}{\bullet}$  &\ftn $\mathbf{+1}~\textcolor{red}{\bullet}$ \\ \cline{2-7} \cline{9-13}
        & \multicolumn{1}{|c|}{\ftn$\Phi_{H,FB}$}& \ftn$-22~\textcolor{red}{\bullet}$ & \ftn$-10~\textcolor{blue}{\bullet}$& \ftn$-1~\textcolor{red}{\bullet}$ & \ftn$+1~\textcolor{red}{\bullet}$& \ftn$-1~\textcolor{red}{\bullet}$ & & \ftn$-42~\textcolor{red}{\bullet}$ & \ftn$-31~\textcolor{blue}{\bullet}$& \ftn$-16~\textcolor{blue}{\bullet}$ & \ftn$+5~\textcolor{red}{\bullet}$& \ftn$+2~\textcolor{red}{\bullet}$\\ \cline{2-7} \cline{9-13}
        & \multicolumn{1}{|c|}{\ftn$\Phi_{H,FISTA}$}& \ftn$-21~\textcolor{red}{\bullet}$&\ftn$\mathbf{-12}~\textcolor{blue}{\bullet}$& \ftn$\mathbf{-10}~\textcolor{blue}{\bullet}$&\ftn$-1~\textcolor{blue}{\bullet}$& \ftn$\mathbf{-2}~\textcolor{blue}{\bullet}$   &  & \ftn$\mathbf{-42}~\textcolor{red}{\bullet}$&\ftn$\mathbf{-31}~\textcolor{blue}{\bullet}$& \ftn$\mathbf{-22}~\textcolor{blue}{\bullet}$&\ftn$+7~\textcolor{red}{\bullet}$& \ftn$+2~\textcolor{red}{\bullet}$  \\ \hline
        \end{tabular}
\caption[Caption for LOF]{\label{tab:blur_comparison} \small For each degradation : the first line of each subtable represents the computation time (in sec) needed by FISTA to reach $5,2,1,0.1$ and $0.01\%$ of the distance $\Vert F_h(x_{h,0})-F_h(x_{h,*})\Vert$. Then for each type of minimization algorithm at coarse level, we display the CPU time relative  to FISTA \eqref{time} (in $\%$) for the best configuration with a colored bullet : $p=1$ \textcolor{red}{$\bullet$} and $p=2$ \textcolor{blue}{$\bullet$}. In all cases : $m=5$. SNR of $z$ : (1a) $11.05$ (1b) $9.64$ (2a) $11.03$ (2b) $9.63$. SNR of $x_{h,300}$ computed by MMFISTA : (1a) $12.71$ (1b) $11.02$ (2a) $12$ (2b) $10.6$. }
\vspace{-1.5em}
\end{table*}

\section{Convergence of the iterates}\label{sec:convergence}
Provided that we use the coarse models a finite number of times, we can prove the convergence of the iterates to a minimizer of $F=F_h$ and that the rate of convergence remains $O (1/k^2)$. First, we consider the sequence of corrections from the coarse models.
\vspace{-1em}
\begin{lemma} \label{lemma_correction} 
Let $L_{f,h}$ and $L_{f,H}$ the Lipschitz constants of $f_h$ and $f_H$, respectively. Let $\tauhk, \tau_{H,k,\ell}\in (0, +\infty)$ the step sizes taken at fine and coarse levels, respectively.
Assume that $\sup_{k,\ell} \tau_{H,k,\ell} < (L_{f,H})^{-1}$ and that $\sup_{k\in \mathbb{N}}\tauhk < L_{f_h}^{-1}$ and denote \np{$\widehat{\tau}_h = \sup_{k} \taubarh$.}
The sequence \np{$(\ehk)_{k\in \mathbb{N}}$} in $\mathcal{H}$ generated by Algorithm \ref{alg:MMFISTA} defined by : 
\begin{equation*}
   \ehk = 
    (\tauhk)^{-1}\widehat{\tau}_h(\mathrm{Id}-\tauhk\nabla f_h)  \left(I_H^h (\xHm-\xHO) \right) 
    \end{equation*}
if a coarse correction is used at iteration $k$ and $ \ehk = 0 $ otherwise, %
is such that  $\sum_{k\in \mathbb{N}} k\Vert \ehk \Vert < +\infty$.
\end{lemma}
\noindent The proof of this lemma is based on the  fact that if the number of coarse corrections is finite, we only need to construct bounded sequences at coarse level so that $I_H^h(\xHm-\xHO)$ is also bounded.

%
%
\noindent From this result we deduce the following theorem :
\begin{theorem}
\label{th:inertial}
    Consider Algorithm \ref{alg:MMFISTA}, suppose that for all $k \in \mathbb{N}^*$, \np{$t_{h,k}$} in \eqref{eq:fista_steps1}, \eqref{eq:fista_steps2} satisfies (AD) conditions \cite{aujol2015}. Suppose that the assumptions of Lemma \ref{lemma_correction} hold. Then : 
    \vspace{-0.5em}
\begin{itemize}
    \item The sequence $(k^{2}\left(F_h(\xhk)-F_h(x^*)\right))_{k \in \mathbb{N}}$ belongs to $\ell_{\infty}(\mathbb{N})$.
    \vspace{-0.5em}
    \item The sequence $(x_{h,k})_{k \in \mathbb{N}}$ given by Algorithm \ref{alg:MMFISTA} weakly converges to a minimizer of $F_h$.
\end{itemize}
\vspace{-1em}
\begin{proof}
    We combine \cite[Theorem 3.5, 4.1, and Corollary 3.8]{aujol2015} with Lemma  \ref{lemma_correction} to prove the desired result.
\end{proof}
\end{theorem}
\vspace{-2em}
\section{Results}
\label{sec:results}

We numerically illustrate the performance of our algorithm in the context of image restoration.%

\begin{figure}[H]
\vspace{-1em}
 \begin{minipage}[c][1\width]{
	   0.5\textwidth}
        \setlength\tabcolsep{2pt}
        \begin{tabular}{ccccc}
        \multirow{4}{*}[-0.3in]{\includegraphics[width=0.27\textwidth]{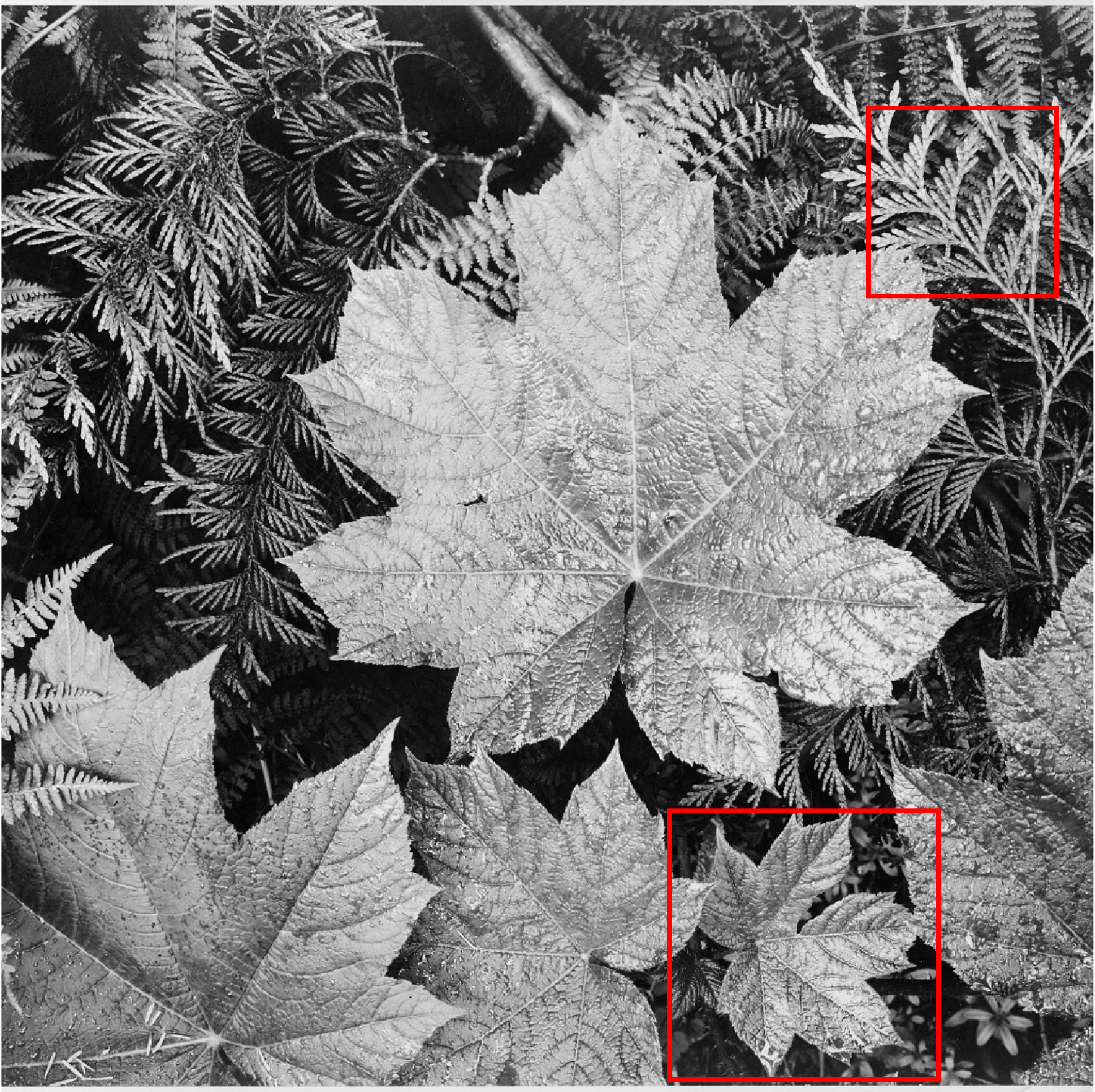}} 
        & \ftn $z$ & \ftn $x^{\mbox{\tiny FISTA}}_{h,2}$ & \ftn $z$ & \ftn $x^{\mbox{\tiny FISTA}}_{h,2}$\\
        & \includegraphics[height=0.14\textwidth]{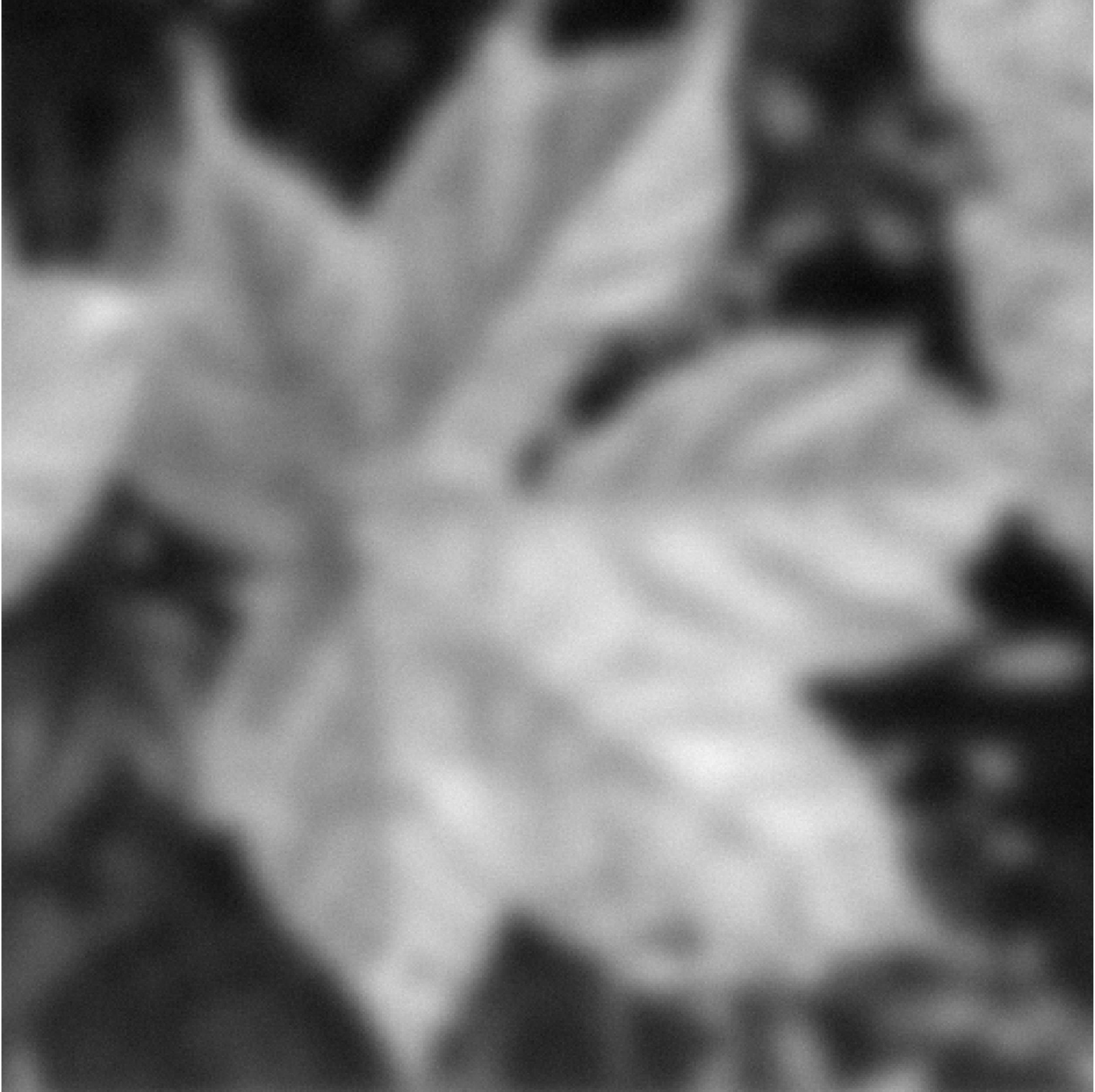}   & \includegraphics[height=0.14\textwidth]{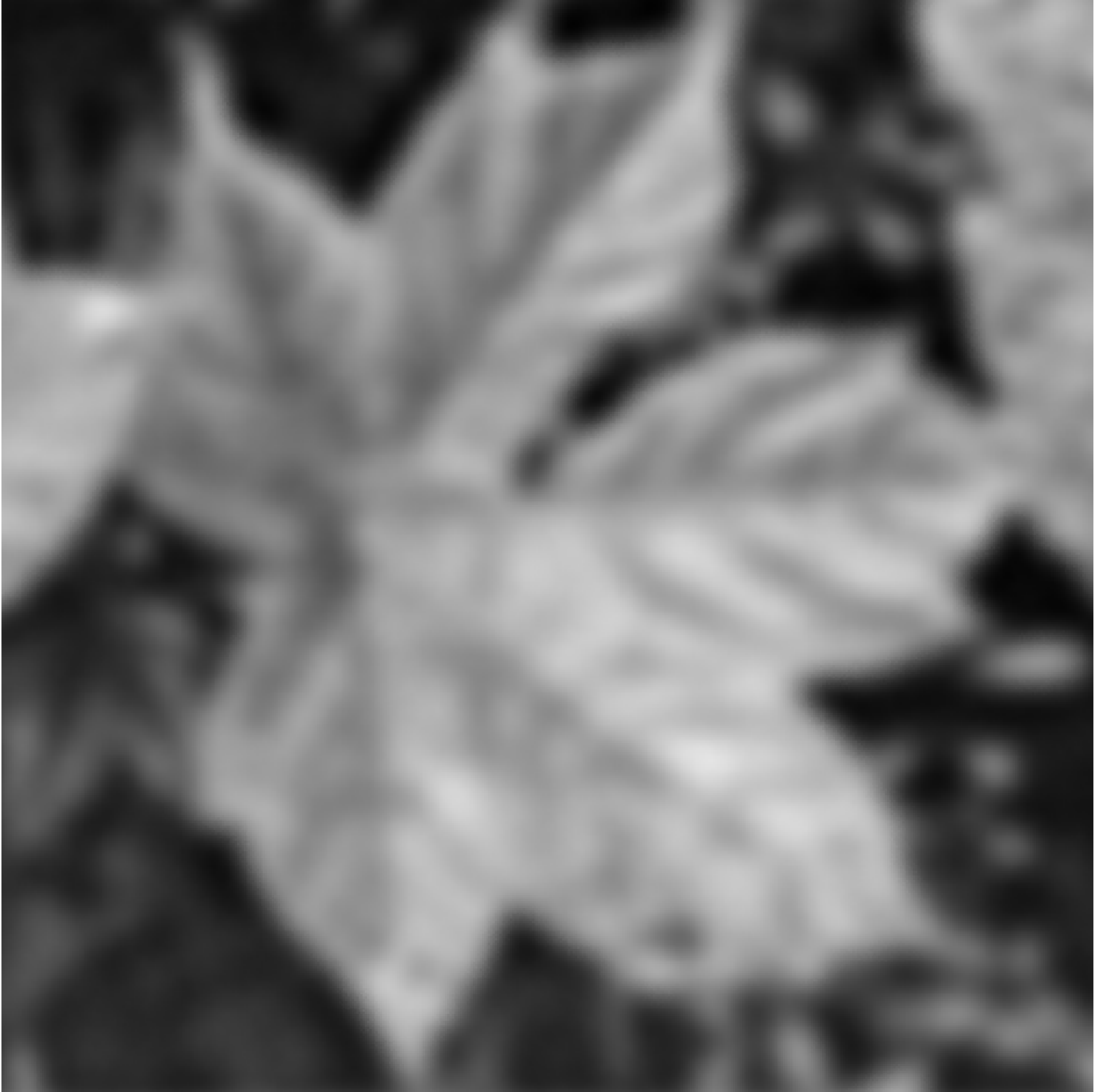}  & \includegraphics[width=0.14\textwidth]{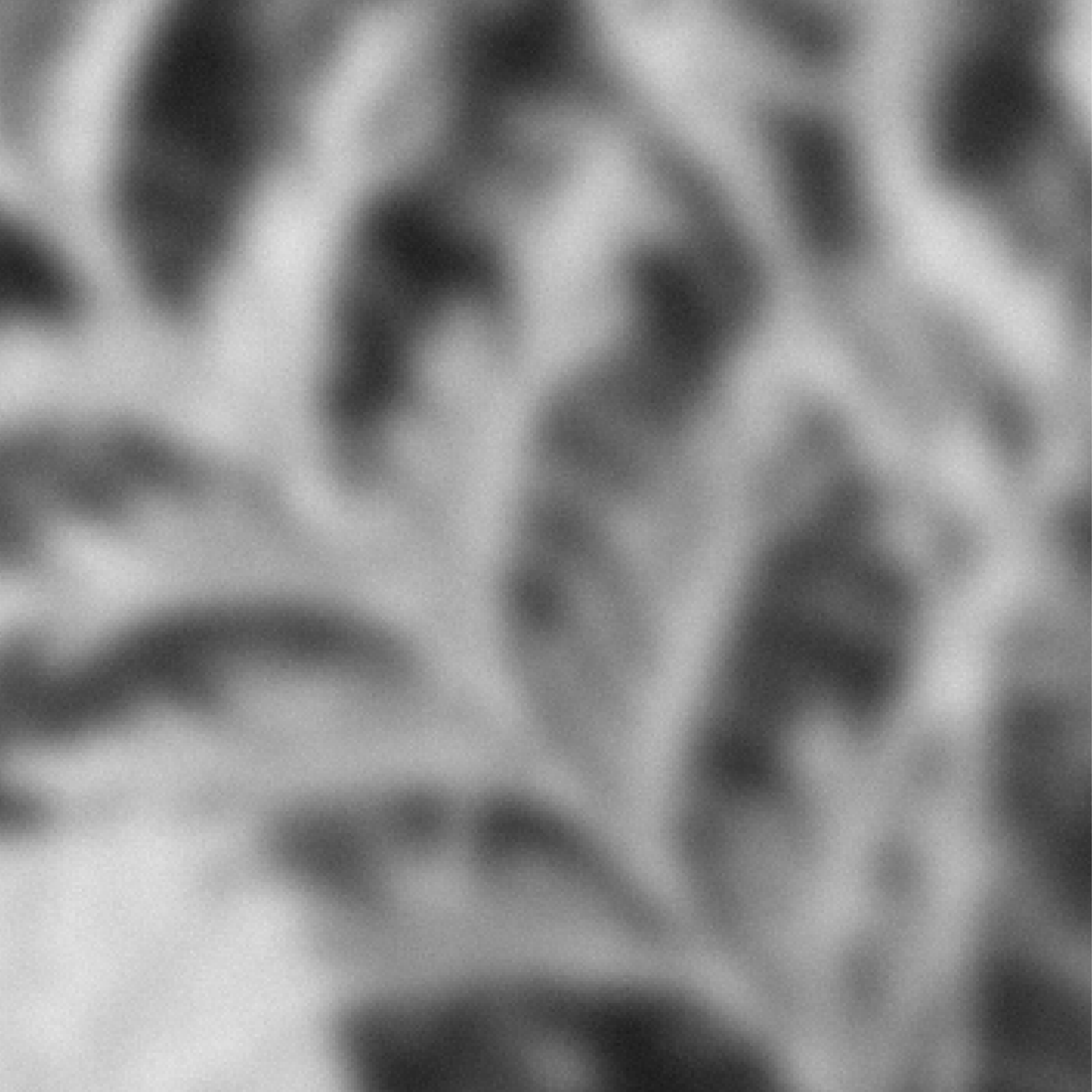}  & \includegraphics[height=0.14\textwidth]{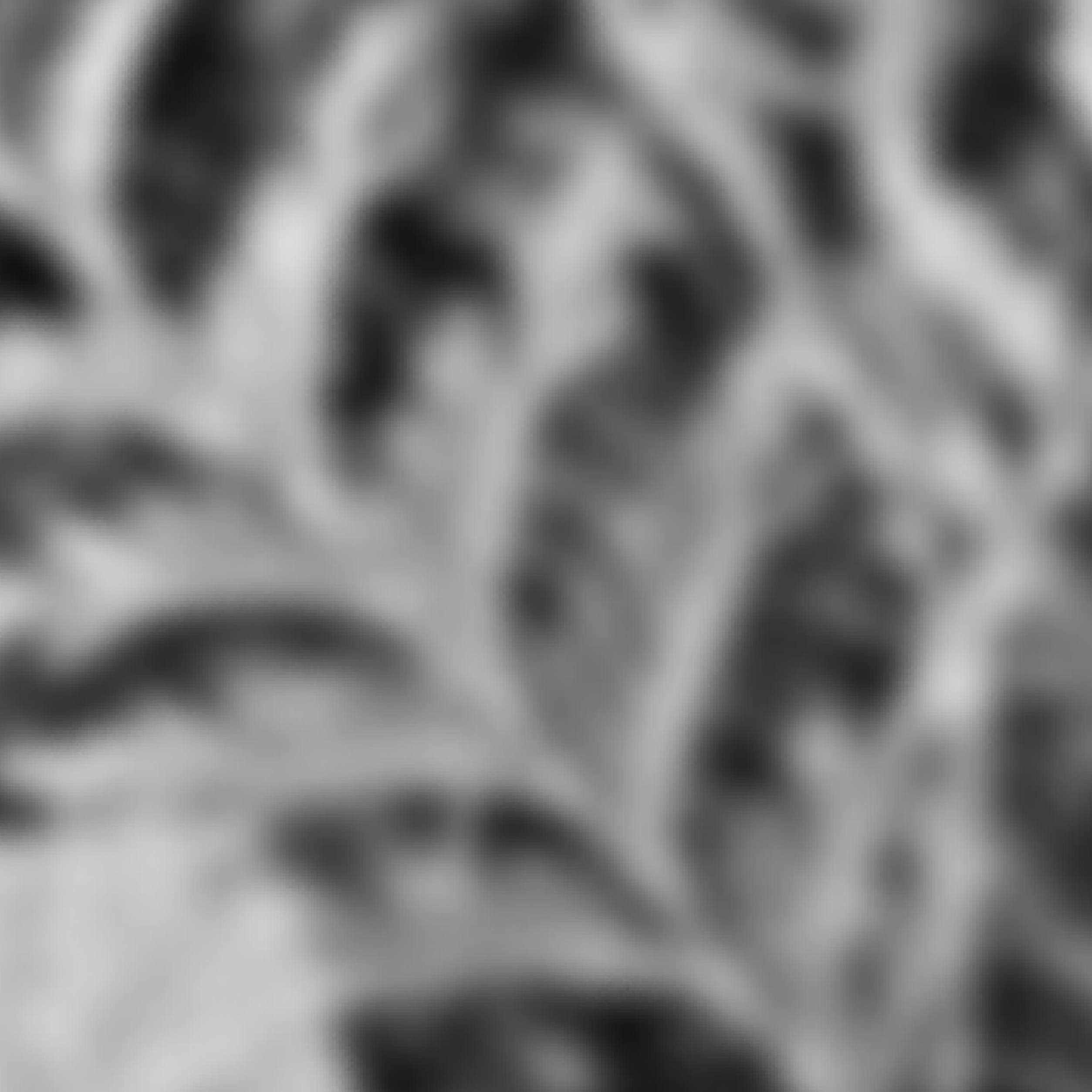}  \\
        & \ftn $x^{\mbox{\tiny MMFISTA}}_{h,300}$ & \ftn $x^{\mbox{\tiny MMFISTA}}_{h,2}$ & \ftn $x^{\mbox{\tiny MMFISTA}}_{h,300}$ & \ftn  $x^{\mbox{\tiny MMFISTA}}_{h,2}$\\
        & \includegraphics[height=0.14\textwidth]{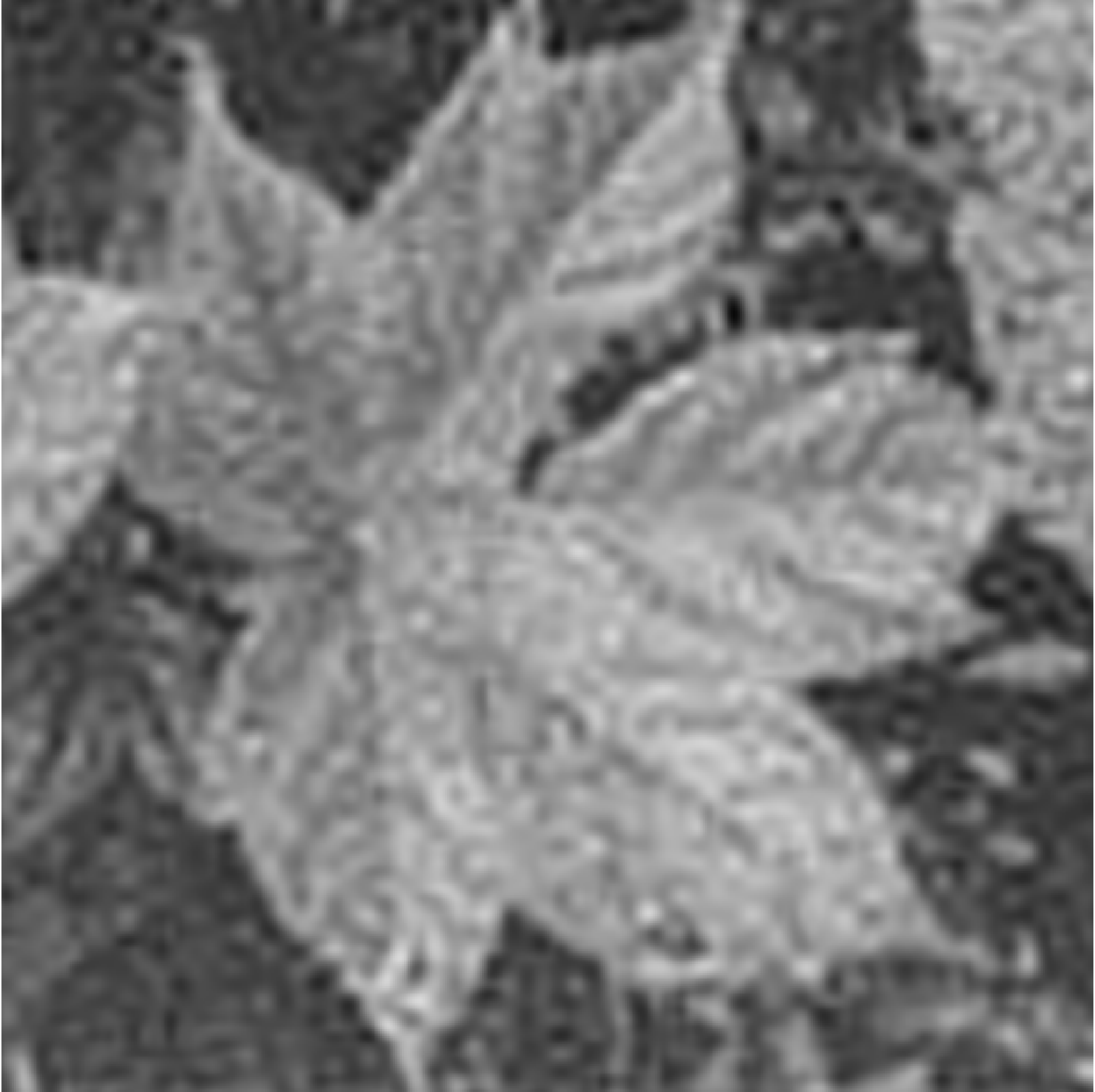}  & \includegraphics[height=0.14\textwidth]{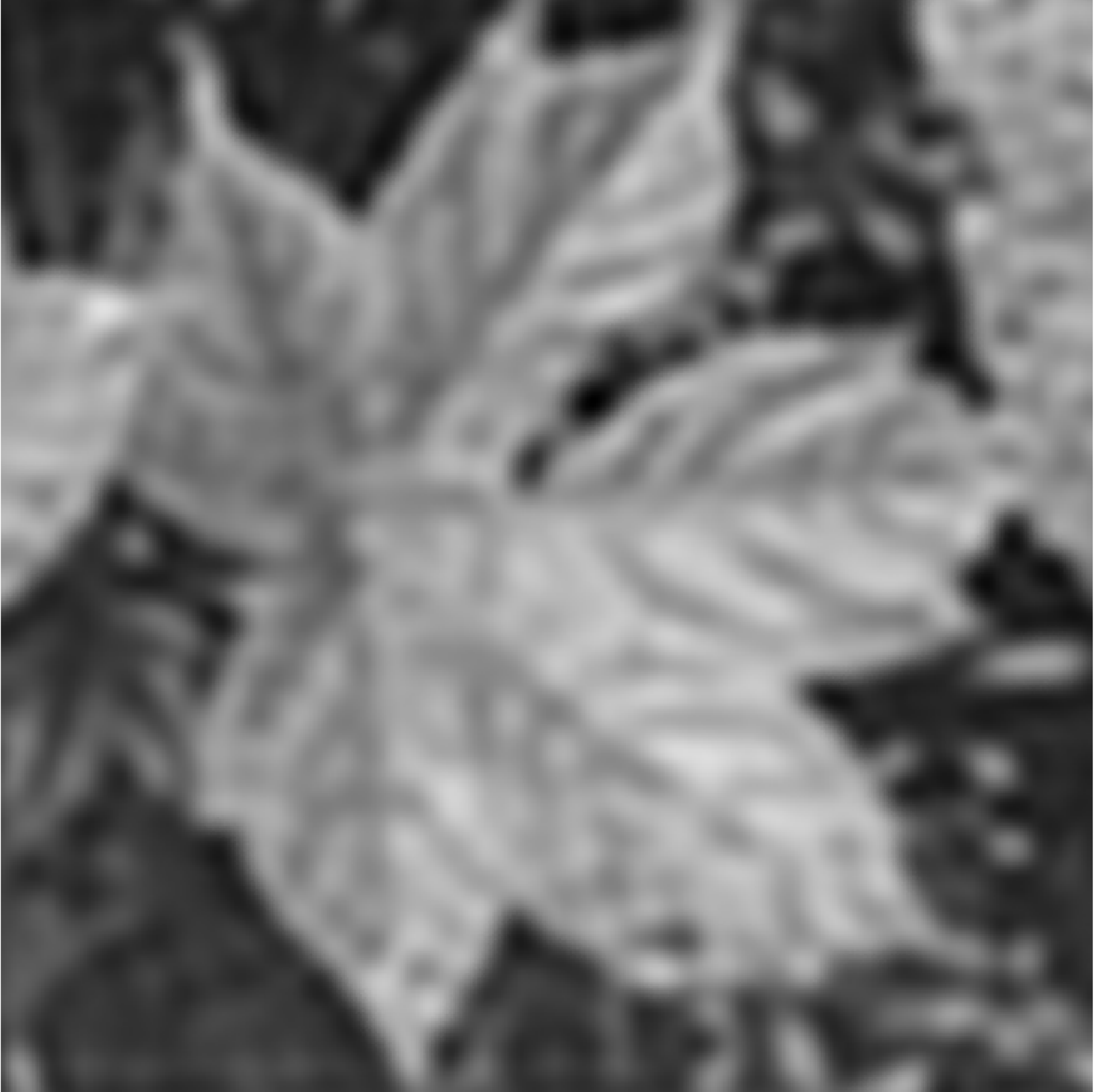} & \includegraphics[width=0.14\textwidth]{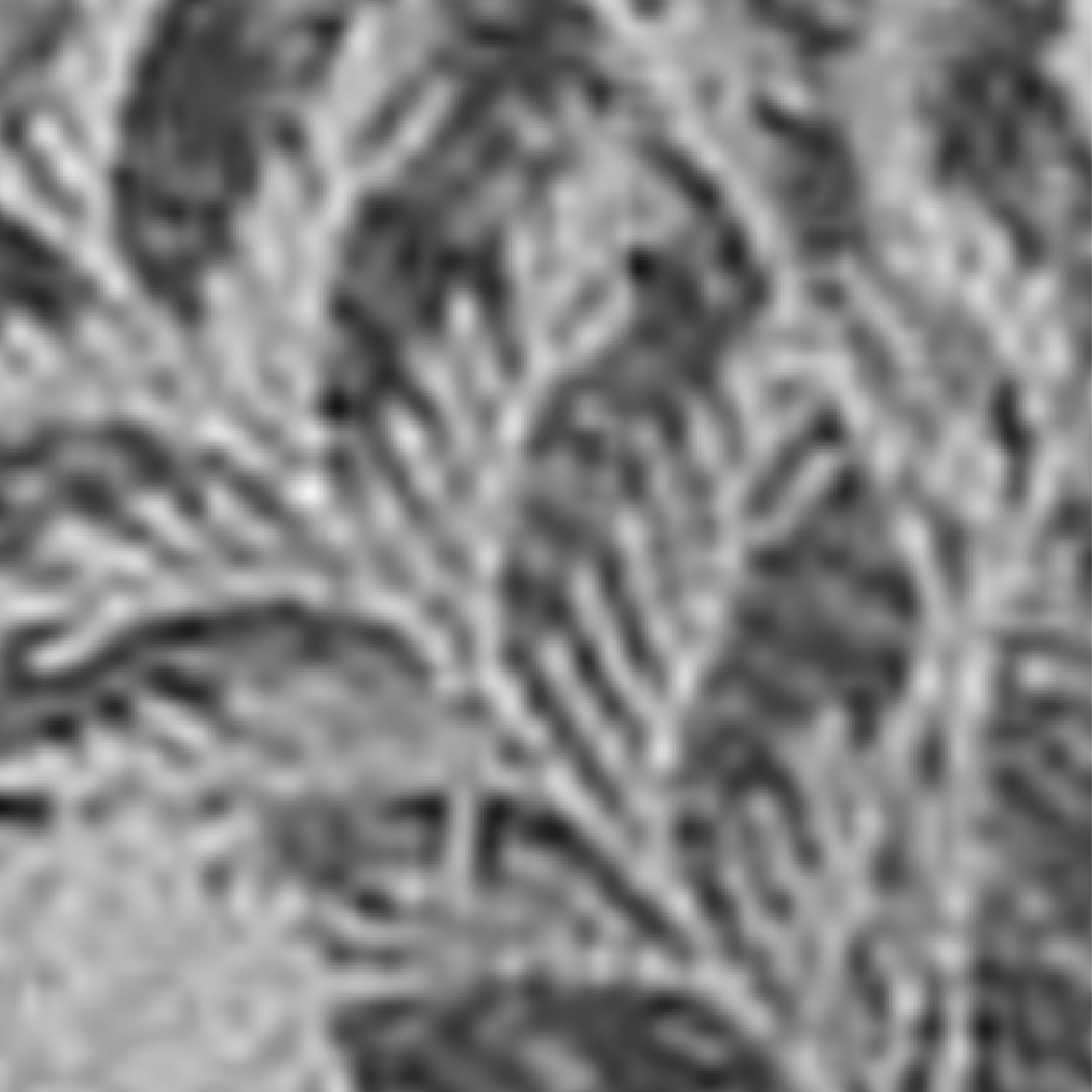} & \includegraphics[height=0.14\textwidth]{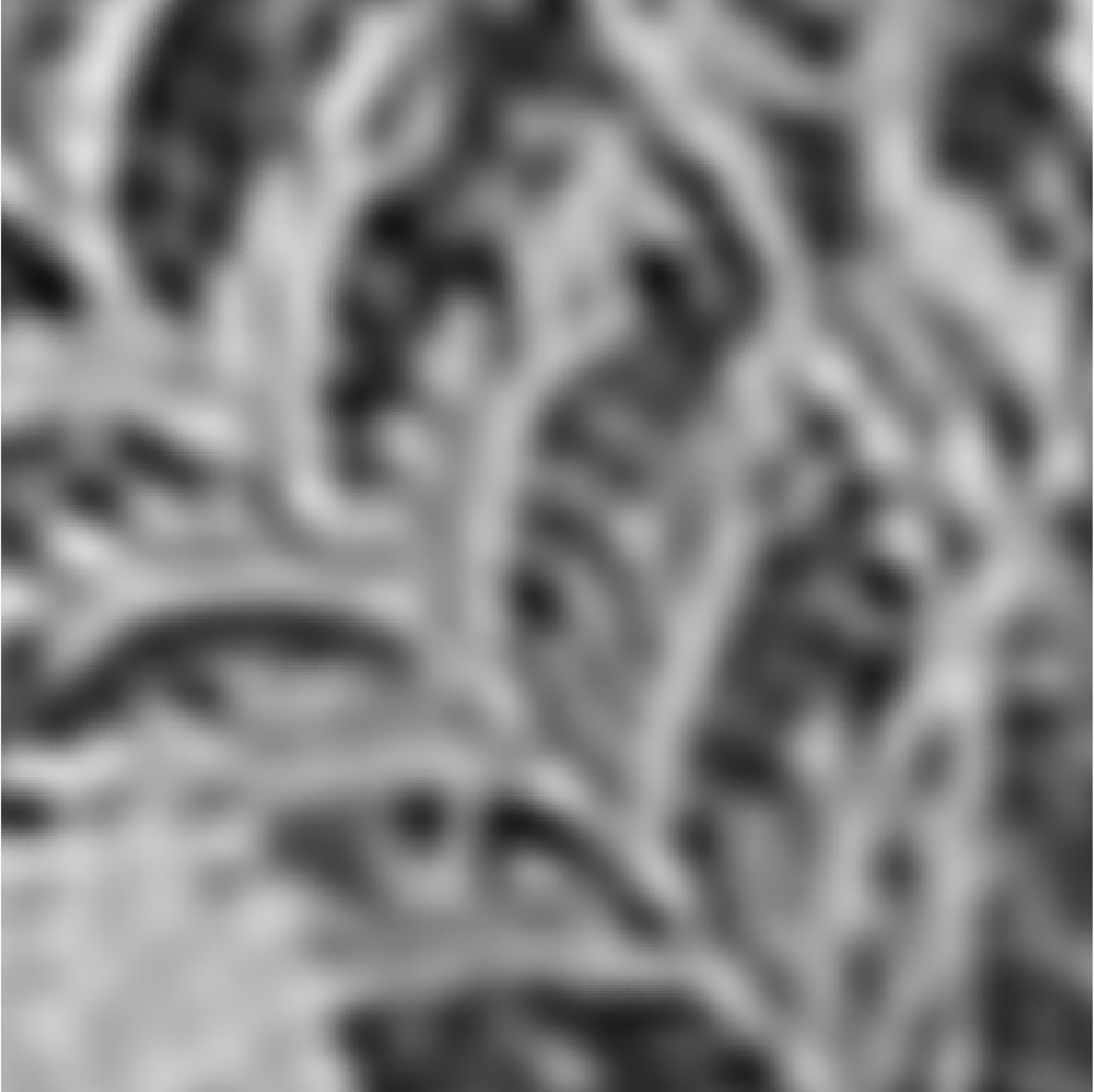}
        \end{tabular}
        \centering
        \setlength\tabcolsep{1pt}
        \begin{tabular}{ll}
        \centering
        \includegraphics[width=0.45\textwidth]{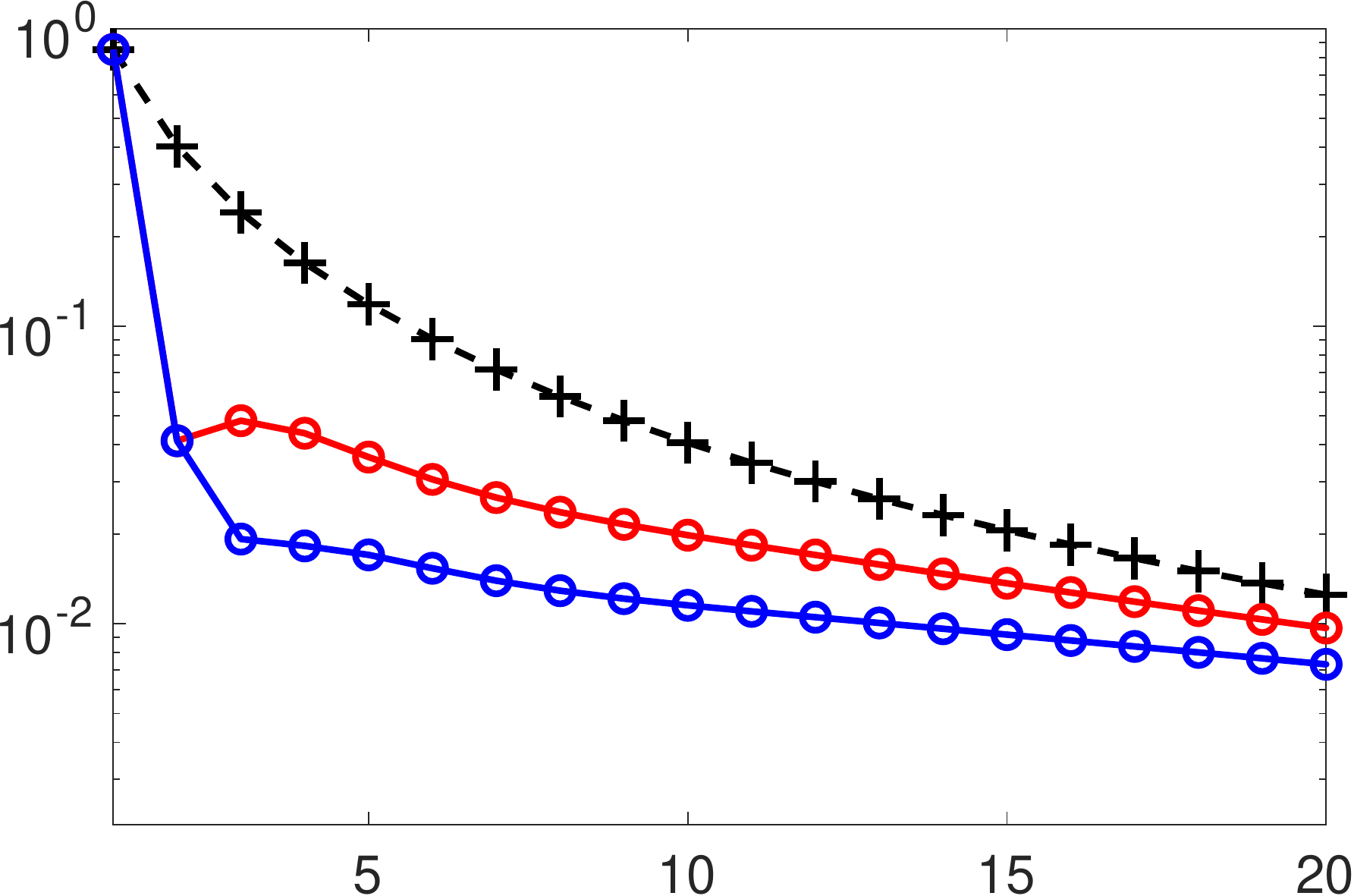}    & \includegraphics[width=0.45\textwidth]{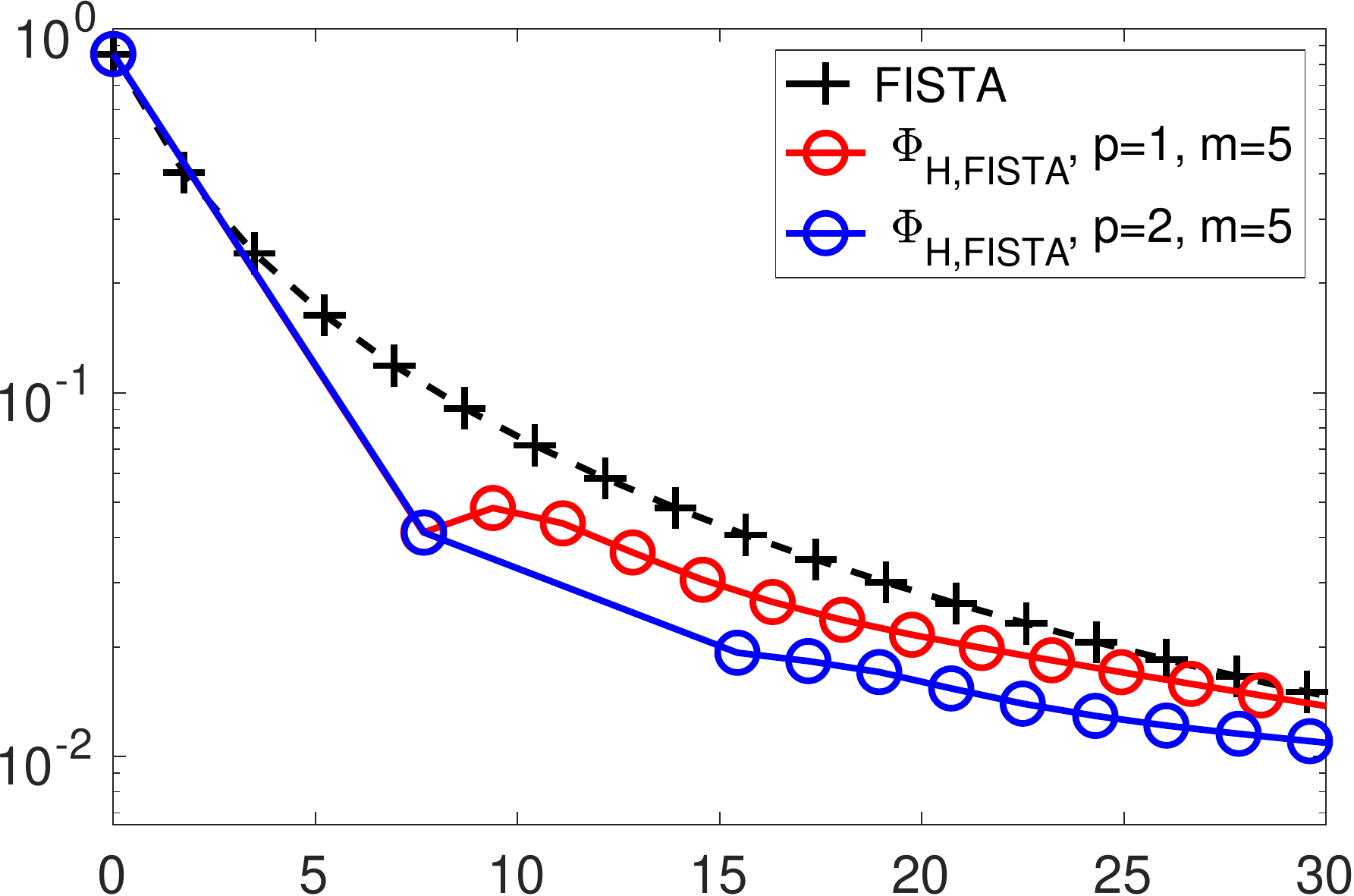}
        \end{tabular}
\end{minipage}

\vspace{-3.8em}
\caption[Caption for LOF]{\label{fig:ansel} \footnotesize{\textit{Top} : From left to right : Original $2048 \times 2048$ image\footnotemark \np{$\overline{x}$ , (first row) zoom} of the degraded image $z$ %
\np{for a noise with $\sigma = 0.01$ and a Gaussian blur of size $40\times 40$ and 7.3 standard deviation and of $x_{h,2}$ computed by FISTA. (second row)} zoom of $x_{h,2}$ and $x_{h,300}$ computed by MMFISTA. \textit{Bottom}: (left) Evolution of $F_h$ versus iterations for MMFISTA with $\Phi_{H,FISTA}$ for $p=1,2$, $m=5$ ; (right) Same for CPU time (in sec). $\lambda_h = 1.7\times 10^{-4}$.}}
\end{figure}

\noindent\textbf{Dataset and degradation -- }We consider large images (Fig.\ref{fig:ansel}) : $2048 \times 2048$, i.e., $N = (2^{J})^2$ with $J=11$, giving  $N \simeq 4\times10^6$. The linear degradation operator $A_h$ is constructed with HNO \cite{hansen2006} as a Kroenecker product with Neumann boundary conditions and we add a Gaussian noise (see the legend of Fig.\ref{fig:ansel} for details). In all tests, the regularization parameter $\lambda_h$ was chosen by a grid search, in order to maximize the SNR of $\widehat{x}$ computed by FISTA at convergence. Also, we choose $x_0$ as the Wiener filtering of $z$.

\noindent\textbf{Multilevel architecture -- }We use a 5-levels hierarchy: from $2048 \times 2048$ ($J=11$) to $128 \times 128$ (indexed by $J-4$). %
We choose $I_h^H$ as the low scale projection on a symlet wavelet with $10$ vanishing moments and $I_H^h = \frac{1}{4} (I_h^H)^T$. We then construct $f_H$ with the blurring matrix $A_H = I_h^H A_h I_H^h$ (which is never used explicitly due to the properties of the Kroenecker product \cite{hansen2006, parpas2017}). Thus $f_h = \frac{1}{2}\Vert A_h x_h -z_h\Vert^2$ and $f_H = \frac{1}{2}\Vert A_H x_H - I_h^H z_h\Vert^2$. The penalty term $g_h = \Vert W_h x_h \Vert_1$ is defined with a full wavelet decomposition over $J$ le\-vels, we construct $g_H = \Vert W_H x_H \Vert_1$ with a decomposition over $J-1$ levels, until $J-4$ levels with $\lambda_H = \lambda_h/4$. The Moreau envelope parameter associated with $g_H$ is set to %
$\gamma_H = 1.1$ while $\gamma_h$ is set to $1$, but both values do not seem to be critical here. 
\noindent 

\noindent \textbf{Visual result --} We display the restored image $\widehat{x}$ and the convergence curves as a function of the iterations and the CPU time for one case in Fig.\ref{fig:ansel}. For clarity, we only display the behaviour of the method with $\Phi_{H,FISTA}$.

\noindent\textbf{Performance assessment -- } We measure $\mbox{Time}_{\mbox{\footnotesize{MMFISTA}}}$, the CPU time needed to reach $5,2,1,0.1$ and $0.01\%$ of the distance %
\np{$\Vert F_h(x_{h,0})-F_h(\widehat{x})\Vert$, }
where %
\np{$\widehat{x}$}
is computed beforehand by FISTA, and we compare it to $\mbox{Time}_{\mbox{\footnotesize{FISTA}}}$, the CPU time of FISTA.  
We tested the performance for several values of $m$, and among our numerous numerical experiments, $m=5$ at the different coarse levels appears to be a good compromise whatever the noise and blur levels.  We report in Tab.\ref{tab:blur_comparison} the  
\begin{equation}\label{time}
\frac{\mbox{Time}_{\mbox{\footnotesize{MMFISTA}}} - \mbox{Time}_{\mbox{\footnotesize{FISTA}}}}{\mbox{Time}_{\mbox{\footnotesize{FISTA}}}}\times 100,
\end{equation}
for $m=5$ at every coarse levels. In this table we evaluate:\\
$\bullet$ The impact of $p$. In our numerical experiments we only consider $p=1$ (\textcolor{red}{$\bullet$}) or $p=2$ (\textcolor{blue}{$\bullet$}) uses of the coarse models, performed at the beginning of the iterative process. They allow to quickly determine the low frequencies components of the solution at the fine level. The choice of $p$ depends on the sought accuracy. 
If a rough approximation is sufficient, fixing $p=1$ is the best choice, while $p=2$ is better for lower thresholds. While we obtain good gains for those, for very low ones the use of a multilevel strategy is not useful, but note that it doesn't deteriorate the performance either. \\
$\bullet$ The impact of noise and blur level. 
For all methods acceleration increases significantly as the blur gets worse. Moreover, as the noise decreases, the improvement obtained with $\Phi_{H,FISTA}$ as compared to others $\Phi_H$ increases. 

The main takeaway from these experiments is that with a few coarse corrections, our method can provide good approximations of the solution while staying competitive with FISTA for high precision approximations.
\vspace{-3em}
\footnotetext{A close-up of leaves in Glacier National Park, Montana taken by Ansel Adams in the 1930s}
\vspace{-2.2em}

\section{Conclusion}
\vspace{-0.3em}
We have proposed a convergent multilevel FISTA method for image restoration that 
reaches rough approximations of the optimal solution in a much smaller CPU time than FISTA and to consider large images. A future research perspective is to extend this approach to other proximal algorithmic schemes and to study/or improve the associated convergence rates. We also want to investigate the influence of the information transfer operators, which remains an open question.
\newpage
\bibliographystyle{IEEEbib}
\bibliography{strings,refs}
\end{document}